\documentclass[a4paper,10pt]{article}
\usepackage[DIV=10]{typearea}
\usepackage{hyperref}
\usepackage{microtype}
\usepackage{graphicx}
\usepackage{amssymb}
\usepackage{amsmath}

\newtheorem{theorem}{Theorem}[section]
\newtheorem{proposition}{Proposition}[section]

\newtheorem{e-proposition}[theorem]{Proposition}

\newtheorem{definition}[theorem]{Definition}

\begin{document}
\title{The method of freezing as a new tool for nonlinear reduced basis approximation of parameterized evolution equations}
\author{Mario Ohlberger\footnote{Institute for Computational and Applied Mathematics \& Center for Nonlinear Science, University of M\"unster, Einsteinstr.\ 62, 48149 M\"unster, Germany, mario.ohlberger@uni-muenster.de}, Stephan Rave\footnote{Institute for Computational and Applied Mathematics, University of M\"unster, Einsteinstr.\ 62, 48149 M\"unster, Germany, stephan.rave@uni-muenster.de}}
\date{April 14, 2013}
\maketitle

\begin{abstract}
We present a new method for the nonlinear approximation of the solution manifolds of parameterized nonlinear evolution problems, in particular in hyperbolic regimes with moving discontinuities.
Given the action of a Lie group on the solution space, the original problem is reformulated as a partial differential algebraic equation system by decomposing the solution into a group component and a spatial shape component, imposing appropriate algebraic constraints on the decomposition.
The system is then projected onto a reduced basis space.
We show that efficient online evaluation of the scheme is possible and study a numerical example showing its strongly improved performance in comparison to a scheme without freezing.
\end{abstract}

\vspace*{\baselineskip}
\vspace*{\baselineskip}


\vspace*{\baselineskip}
\vspace*{\baselineskip}

\section{Introduction}
Reduced basis (RB) methods are effective tools for approximating the solution manifolds of parameterized evolution problems by 
low-dimensional linear spaces, enabling fast online evaluation of the solution for arbitrary parameter values.  
For linear problems, the POD-Greedy algorithm \cite{M2AN08} is by now an established 
reduction approach which has recently been proved to be optimal in the sense that exponential or algebraic 
convergence rates of the Kolmogorov $n$-width are maintained by the algorithm \cite{haasdonk}. 
The approach has been further extended to nonlinear settings in \cite{drohmann} based on empirical interpolation of the nonlinear operators.

For convection dominated problems with low regularity however, the performance of RB-methods is limited by the fact that translation of functions is nonlinear in the translation vector:
the solution has to be approximated at every location in space it is being transported to, resulting in merely linear decay of the $n$-width.
Our aim is therefore to extend the above approaches by additionally allowing transformations of the reduced spaces given
by an appropriate group action on the solution space.

Originally developed for the study of relative equilibria of evolution equations \cite{beyn,rowley}, 
the \emph{method of freezing} allows us to obtain such a decomposition of the solution into a group and shape component 
for arbitrary Lie group actions, provided they satisfy the equivariance condition (\ref{equivariance}).
Combining this method with RB-techniques, we obtain a new nonlinear reduction method for parameterized nonlinear 
evolution equations.

\section{The method of freezing}
Assume we are given a parameter dependent nonlinear Cauchy problem of the form
\begin{equation}\label{problem}
\partial_t u_\mu(t) + \mathcal{L}_\mu(u_\mu(t)) = 0, \qquad
    u_\mu(0) = u_{0}
\end{equation}
where $u_\mu(t) \in V$ is a function of space for each $t \in [0, T]$ and $\mathcal{L}_\mu$ is a partial differential operator acting on $V$. 

Given a Lie group $G$ acting smoothly on $V$ by linear operators, we want to decompose the solution $u_\mu(t)$ into a group component $g_\mu(t) \in G$ and a shape component $v_\mu(t) \in V$ such that
\begin{equation} \label{decomp}
	u_\mu(t) = g_\mu(t)\,.\,v_\mu(t)\,.
\end{equation}
Here the action of $G$ on $V$ is denoted by a lower dot.

As a guiding example consider the action of the Lie groups $G= \mathbb{R}^d$ on functions $v: \mathbb{R}^d \longrightarrow \mathbb{R}$ via translation, i.e. $(g\,.\,v)(x)=v(x - g)$ for $g \in \mathbb{R}^d$. 
In this case, if $u_\mu(t)$ is a solution of (\ref{problem}) drifting along a certain trajectory in space, we want to find a decomposition (\ref{decomp}) such that the drift is captured by the evolution of $g_\mu(t)$ whereas $v_\mu(t)$ becomes as stationary as possible, only describing the change of the shape of $u_\mu(t)$ over time.

Note that, in general, the group action does not have to be induced by a mapping of the underlying spatial domain (e.g. shifts or rotations) but can also involve more general transformations of $V$.

Inserting (\ref{decomp}) into (\ref{problem}), we formally get
\[ \partial_t g_\mu(t)\,.\,v_\mu(t) + g_\mu(t)\,.\,\partial_t v_\mu(t) + \mathcal{L}_\mu(g_\mu(t)\,.\,v_\mu(t)) = 0 \]
which can be rewritten as
\begin{equation}
\partial_t v_\mu(t) + g_\mu(t)^{-1}\,.\,\mathcal{L}_\mu(g_\mu(t)\,.\,v_\mu(t)) + \mathfrak{g}_\mu(t)\,.\,v_\mu(t) = 0,
  \qquad
	  \mathfrak{g}_\mu(t) = g_\mu(t)^{-1}\cdot \partial_t g_\mu(t)
\end{equation}
where $\mathfrak{g}_\mu(t)$ is an element of the Lie algebra $\mathrm{L}G$ of $G$, i.e. the tangential space of $G$ at its neutral element $1_G$.

Since the decomposition ($\ref{decomp}$) introduces $\dim G$ additional degrees of freedom, the resulting system is now underdetermined. 
It is the main idea of the method of freezing to compensate for these degrees of freedom by adding appropriate algebraic constraints $\Phi$ which force $v_\mu$ to have minimal change over time. 
These constraints are called \emph{phase conditions}.
Thus, the group component of the solution is automatically determined by the phase condition, and no a-priori knowledge of the evolution of the solution is necessary. 

If one further assumes that the operator $\mathcal{L}_\mu$ is \emph{equivariant} under the group action, i.e.
\begin{equation}\label{equivariance}
g^{-1}\,.\,\mathcal{L}_\mu(g\,.\,v) = \mathcal{L}_\mu(v) \qquad \forall g \in G, v \in V,
\end{equation}
the system decuples into the partial differential algebraic equation (PDAE) of index one
\begin{equation}\label{frozenproblem}
	  \partial_t v_\mu(t) + \mathcal{L}^G_{\mu, \mathfrak{g}_\mu(t)}(v_\mu(t)) = 0,
	  \qquad
	  \Phi(v_\mu(t), \mathfrak{g}_\mu(t)) = 0
\end{equation}
with $\mathcal{L}^G_{\mu, \mathfrak{h}}(u)=\mathcal{L}_\mu(u) + \mathfrak{h}\,.\,u$ and the ordinary differential equation
  $\partial_t g_\mu(t) = g_\mu(t)\cdot \mathfrak{g}_\mu(t)$,
called the \emph{reconstruction equation}.
The initial conditions are $v_\mu(0) = u_{0}$ and $g_\mu(0) = 1_G$.

Different choices of phase conditions are possible \cite{beyn}. 
We will restrict ourselves here to the so-called \emph{orthogonality phase condition}:
assuming that $V$ is equipped with an inner product, we require that the evolution of $v_\mu$ is at each point in time orthogonal to the action of $\mathrm{L}G$, i.e. $(\partial_t v_\mu(t),\, \mathfrak{h}\,.\,v_\mu(t)) = 0$ for all $\mathfrak{h} \in \mathrm{L}G$. 
Inserting (\ref{frozenproblem}), we obtain 
\[ \bigl(\mathcal{L}_\mu(v_\mu(t)),\, \mathfrak{h}\,.\,v_\mu(t)\bigr) + \bigl(\mathfrak{g}_\mu(t)\,.\,v_\mu(t),\, \mathfrak{h}\,.\,v_\mu(t)\bigr)  = 0\qquad \forall\mathfrak{h} \in \mathrm{L}G. \]
After choosing a basis for $\mathrm{L}G$, this leads to a linear $\dim G \times \dim G$ equation system for $\mathfrak{g}_\mu(t)$.
Denoting the basis vectors by $\mathfrak{e}_1,\ldots,\mathfrak{e}_{\dim G}$, this system can be written as
\begin{equation*}
	\Phi(v_\mu(t), \mathfrak{g}_\mu(t))= \Bigl[ \bigl(\mathfrak{e}_r\,.\,v_\mu(t),\, \mathfrak{e}_s\,.\,v_\mu(t)\bigr) \Bigr]_{r,s} \cdot \Bigl[ \mathfrak{g}_{\mu,s}(t) \Bigr]_s
	+ \Bigl[ \bigl(\mathcal{L}_\mu(v_\mu(t)),\, \mathfrak{e}_r\,.\,v_\mu(t)\bigr) \Bigr]_r  = 0
\end{equation*}
where $\mathfrak{g}_{\mu,s}$ denotes the $s$th component of $\mathfrak{g}_{\mu}$ with respect to the chosen basis.

Depending on the specific problem, many different choices of space-time discretizations of the frozen system (\ref{frozenproblem}) are possible.
To keep the notation simple, we will restrict ourselves to an explicit Euler time-stepping:

\begin{definition}\label{def:frozenscheme}
Let an $H$-dimensional space $V_H$ and discrete parameter dependent operators $\mathbb{L}_\mu$, $\mathbb{L}^G_{\mu,\mathfrak{h}}$ and $\mathbb{S}^G_r$ on $V_H$ be given, approximating the operators $\mathcal{L}_\mu$, $\mathcal{L}_{\mu,\mathfrak{h}}$ and $\mathfrak{e}_r\,.\,(\,\cdot\,)$. 
Since $\mathrm{L}G$ acts on $V$ by linear operators, we can assume that the $\mathbb{S}^G_i$ are linear.
Furthermore let $P_H:V \longrightarrow V_H$ be a projection operator.
For $K \in \mathbb{N}$, $\Delta t = T/K$ and $0 \leq k \leq K$, we define the discrete solutions $v^k_\mu \in V_H$, $\mathfrak{g}_\mu^k \in \mathrm{L}G$ of the frozen system (\ref{frozenproblem}) by the equations
$v_\mu^0 = P_H(u_0)$
and
\begin{equation}\label{discretescheme}
     v_\mu^{k+1} = v_\mu^k - \Delta t\, \mathbb{L}^G_{\mu,\mathfrak{g}^k_\mu}(v_\mu^k),
     \qquad
    \Bigl[ \bigl(\mathbb{S}^G_r(v_\mu^k),\, \mathbb{S}^G_s(v_\mu^k)\bigr) \Bigr]_{r,s} \cdot \Bigl[ \mathfrak{g}_{\mu,s}^k \Bigr]_s
     = 
     - \Bigl[ \bigl(\mathbb{L}_\mu(v_\mu^k),\, \mathbb{S}_r(v_\mu^k)\bigr) \Bigr]_r 
\end{equation}
with $k = 0,\ldots,K - 1$. Moreover, let
\[ g_\mu^0 = 1_G, \qquad g_\mu^{k+1}= g_\mu^k \cdot \exp_G(\Delta t\cdot \mathfrak{g}_\mu^k) \qquad k = 1,\ldots,n_t - 1. \]
\end{definition}
A discrete approximation of the solution $u_\mu(t)$ is then obtained through
$ u^k_\mu = g^k_\mu\,.\,v^k_\mu $,
given an appropriate discrete approximation of the action of $G$.

\section{The FrozenRB-scheme}
After freezing and discretizing the original problem (\ref{problem}), we will now reduce the resulting system (\ref{discretescheme}) by projecting it onto an $N$-dimensional reduced basis space $V_N \subseteq V_H$.
To achieve fast online evaluation of the reduced scheme, we approximate the nonlinear operators $\mathbb{L}_\mu$, $\mathbb{L}_{\mu, \mathfrak{h}}^G$ using the method of \emph{empirical operator interpolation} \cite{drohmann}.
If $\hat{\varphi}_1,\cdots,\hat{\varphi}_H$ denotes a basis of the dual of $V_H$, this method produces indices $q_1,\ldots,q_M$ and vectors $\xi_1,\ldots,\xi_M,\xi^G_1,\ldots,\xi^G_M \in V_H$ such that
\[
	\mathbb{L}_\mu(v) \approx \sum_{m=1}^M \hat{\varphi}_{q_m}(\mathbb{L}_\mu(v))\cdot\xi_m,
	\qquad
	\mathbb{L}^G_{\mu,\mathfrak{h}}(v) \approx \sum_{m=1}^M \hat{\varphi}_{q_m}(\mathbb{L}_{\mu,\mathfrak{h}}^G(v))\cdot\xi_m^G\;.
\]
Given an appropriate projection operator $P_N:V_H \longrightarrow V_N$, the fully reduced scheme then reads
\begin{equation}\label{rbscheme}
	\left\{\begin{array}{l}
	  v_{N,\mu}^{k+1} = v_{N,\mu}^k - \Delta t\, P_N\Bigl(\sum_{m=1}^M\hat{\varphi}_{q_m}(\mathbb{L}^G_{\mu,\mathfrak{g}^k_{N,\mu}}\!\!(v_{N,\mu}^k))\cdot\xi^G_m\Bigr) \\
	  \Bigl[ \bigl(\mathbb{S}^G_r(v_{N,\mu}^k),\, \mathbb{S}^G_s(v_{N,\mu}^k)\bigr) \Bigr]_{r,s} \cdot \Bigl[ \mathfrak{g}_{N,\mu,s}^k \Bigr]_s
     = 
     - \Bigl[ \bigl(\sum_{m=1}^M\hat{\varphi}_{q_m}(\mathbb{L}_\mu(v_{N,\mu}^k))\cdot\xi_m,\, \mathbb{S}_r(v_{N,\mu}^k)\bigr) \Bigr]_r 
  \end{array}\right. 
\end{equation}
with $v^k_{N,\mu} \in V_N$, $\mathfrak{g}^k_{N,\mu}\in\mathrm{L}G$, and $v^0_{N,\mu} = P_N(v^0_\mu)$.

An appropriate reduced space $V_N$ and the interpolation data $q_m$, $\xi_m$, and $\xi^G_m$ can be constructed from solutions of (\ref{discretescheme}) using reduced basis techniques developed in \cite{drohmann}.
Choosing the same interpolation points for $\mathbb{L}_\mu$ and $\mathbb{L}^G_{\mu,\mathfrak{h}}$ ensures that expensive evaluations of nonlinear flux functions have to be carried out only once for both operators.

\subsection{Offline/online decomposition}
The discrete operators $\mathbb{L}_\mu$, $\mathbb{L}_{\mu, \mathfrak{h}}^G$ that arise from standard discretizations share the property of being local in the sense that for an appropriate basis $\varphi_i$ of $V_H$ each evaluation $\hat{\varphi}_i(\mathbb{L}_\mu(v))$ only depends on at most $C$ degrees of freedom of $v$, independently of $i$, $v$ and $H$ ($H$-independent DOF dependence \cite{drohmann}). 
In particular, each of these evaluations can be performed at a speed independent of $H$.

Under this assumption we achieve fast online evaluation of (\ref{rbscheme}) using precomputed data as follows:

\begin{proposition}\label{prop:onlineoffline}
Let $\mathbb{L}_\mu$, $\mathbb{L}_{\mu, \mathfrak{h}}^G$ have the $H$-independent DOF dependence property.
Then there are indices $q^\prime_1,\ldots,q^\prime_L$ with $L \leq 2CM$ such that
\[ 
\hat{\varphi}_{q_m}(\mathbb{L}_\mu(v)) = \hat{\varphi}_{q_m}\Bigl(\mathbb{L}_\mu\Bigl(\sum_{l=1}^L \hat{\varphi}_{q^\prime_l}(v)\cdot\varphi_{q^\prime_l}\Bigr)\Bigr),
\qquad
\hat{\varphi}_{q_m}(\mathbb{L}^G_{\mu,\mathfrak{h}}(v)) = \hat{\varphi}_{q_m}\Bigl(\mathbb{L}^G_{\mu,\mathfrak{h}}\Bigl(\sum_{l=1}^L \hat{\varphi}_{q^\prime_l}(v)\cdot\varphi_{q^\prime_l}\Bigr)\Bigr)
\]
for all $v\in V_H$, $m = 1, \ldots, M$.

If $\psi_1,\ldots,\psi_N$ is a basis of $V_N$, define matrices $\mathbf{P}, \mathbf{EV}, \mathbf{PCL}^{{r,s}}$ and $\mathbf{PCR}^{r}$ for $1 \leq r,s \leq \dim G$ as
\[
\begin{array}{lcl}
	\mathbf{P}_{n,m} &=& \hat{\psi}_n(P_N(\xi^G_m)), \\
	\mathbf{EV}_{l,n} &=& \hat{\varphi}_{q^\prime_l}(\psi_n),
\end{array}
\qquad\quad
\begin{array}{lcl}
	\mathbf{PCL}^{{r,s}}_{n,n^\prime} &=& \bigl(\mathbb{S}^G_r(\varphi_n),\, \mathbb{S}^G_s(\varphi_{n^\prime})\bigr), \\
	\mathbf{PCR}^{r}_{m,n} &=& \bigl(\xi_m,\, \mathbb{S}_r(\varphi_n)\bigr) 
\end{array}
\]
with $1 \leq l \leq L$, $1 \leq m \leq M$ and $1 \leq n, n^\prime \leq N$.
Let moreover $\mathbf{L}_\mu,\mathbf{L}^{G}_{\mu, \mathfrak{g}}: \mathbb{R}^L \longrightarrow \mathbb{R}^M$ be given as
\[ 
[\mathbf{L}_\mu(\mathbf{y})]_m = \hat{\varphi}_{q_m}\Bigl(\mathbb{L}_\mu\Bigl(\sum_{l=1}^L \mathbf{y}_l\cdot\varphi_{q^\prime_l}\Bigr)\Bigr),
\qquad
[\mathbf{L}^{G}_{\mu,\mathfrak{g}}(\mathbf{y})]_m = \hat{\varphi}_{q_m}\Bigl(\mathbb{L}^G_{\mu,\mathfrak{g}}\Bigl(\sum_{l=1}^L \mathbf{y}_l\cdot\varphi_{q^\prime_l}\Bigr)\Bigr).
\]
If $\mathbf{v}^k_{N,\mu}$ is the coefficient vector of $v^k_{N,\mu}$ with respect to the reduced basis of $V_N$, then (\ref{rbscheme}) is equivalent to 
\[
  \left\{\begin{array}{l}
	  \mathbf{v}_{N,\mu}^{k+1} = \mathbf{v}_{N,\mu}^k - \Delta t\, \mathbf{P}\cdot \mathbf{L}^G_{\mu, \mathfrak{g}^k_{N,\mu}}(\mathbf{EV} \cdot \mathbf{v}_{N,\mu}^k)\\
	  \Bigl[ (\mathbf{v}^k_{N,\mu})^T \cdot \mathbf{PCL}^{r,s} \cdot \mathbf{v}^k_{N,\mu}\Bigr]_{r,s} \cdot \Bigl[ \mathfrak{g}_{N,\mu,s}^k \Bigr]_s
     = 
     - \Bigl[ ( \mathbf{L}_{\mu}(\mathbf{EV} \cdot \mathbf{v}_{N,\mu}^k))^T \cdot \mathbf{PCR}^r \cdot \mathbf{v}^k_{N,\mu} \Bigr]_r \ .
  \end{array}\right. 
\]
\end{proposition}

\section{Numerical experiment}

As a first test for our new method we consider the Burgers-type problem from \cite{hyp08}, i.e.\ we solve
\begin{equation}
		\partial_t u_\mu(t) + \nabla \cdot (\mathbf{b}\,u_\mu(t)^\mu) = 0, \qquad
		u_\mu(0) = u_0
	\label{burgers}
\end{equation}
for $t \in [0, 0.3]$ on the domain $\Omega=[0,2]\times [0,1]$ identifying $\{0\}\times [0,1]$ with $\{2\}\times[0,1]$ and $[0,2]\times\{0\}$ with $[0,2]\times\{1\}$.
The parameter $\mu$ is allowed to vary in the interval $[1,2]$.
Moreover, let $\mathbf{b} = [\,1\ 1\,]^T$ and $u_0(x_1,x_2)=1/2(1 + \sin(2\pi x_1)\sin(2\pi x_2))$.

Equation (\ref{burgers}) is invariant under the action of the group $\mathbb{R}^2$ by translations of $\Omega$. 
Since the action of $\mathrm{L}\mathbb{R}^2 = \mathbb{R}^2$ is given by negative spatial derivatives, this leads us to the frozen PDAE
\begin{equation}
\left\{	\begin{array}{l}
		\partial_t v_\mu(t) + \nabla \cdot (\mathbf{b}\,v_\mu(t)^\mu) - \mathfrak{g}_{\mu}(t) \cdot \nabla v_\mu(t) = 0 \\
		\Bigl[ \bigl(\partial_{x_r}v_\mu(t),\, \partial_{x_s}v_\mu(t)\bigr) \Bigr]_{r,s} \cdot \Bigl[ \mathfrak{g}_{\mu,s}(t) \Bigr]_s
		= 
		\Bigl[ \bigl(\nabla \cdot (\mathbf{b}\,v_\mu(t)^\mu),\, \partial_{x_r}v_\mu(t)\bigr) \Bigr]_r \qquad 1\leq r,s \leq 2
	\end{array}
	\label{burgersfrozen}\right.
\end{equation}
with initial condition $v_\mu(0) = u_0$. 
The reconstruction equation for $g_\mu(t)$ is given by $\partial_t g_\mu(t) = \mathfrak{g}_\mu(t)$. 

For the discretization of (\ref{burgersfrozen}) we choose a finite volume scheme on a $120\times 60$ grid and take 100 steps in time.
We use the POD-Greedy and EI-Greedy algorithms from \cite{drohmann} on a fixed snapshot set for the generation of the reduced basis and the operator interpolation data, choosing for simplicity a common interpolation basis $\xi_m = \xi^G_m$ for the operators.

\begin{figure}[p] 
\begin{center} \includegraphics{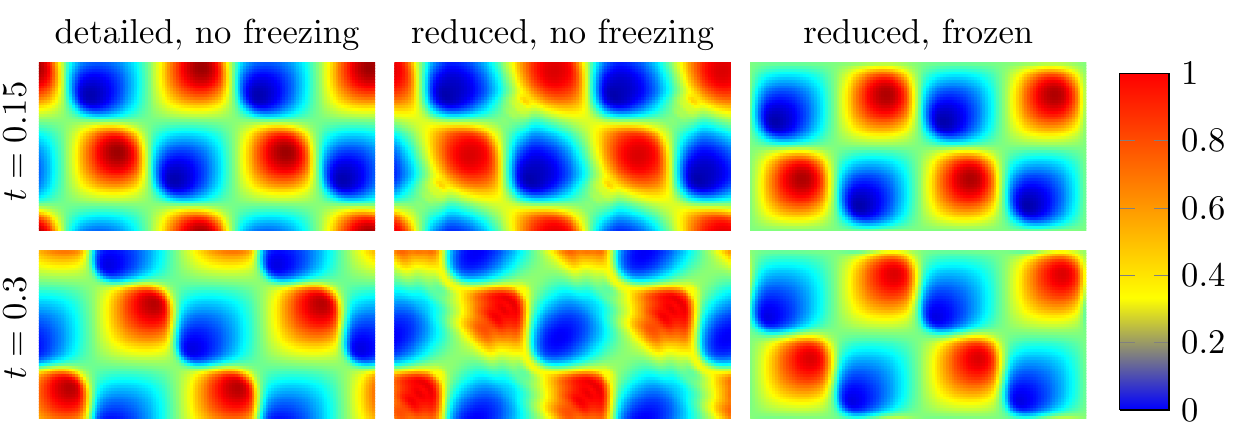} \end{center}
\caption{Numerical solutions of the Burgers problem (\ref{burgers}) for $\mu = 1.5$; left column: detailed approximation; center column: reduced approximation; right column: reduced approximation of the frozen solution $v_\mu$ ($20$ basis vectors, $38$ interpolation points).}
\label{fig:solution}
\end{figure}

\begin{figure}[p] 
\begin{center} \includegraphics{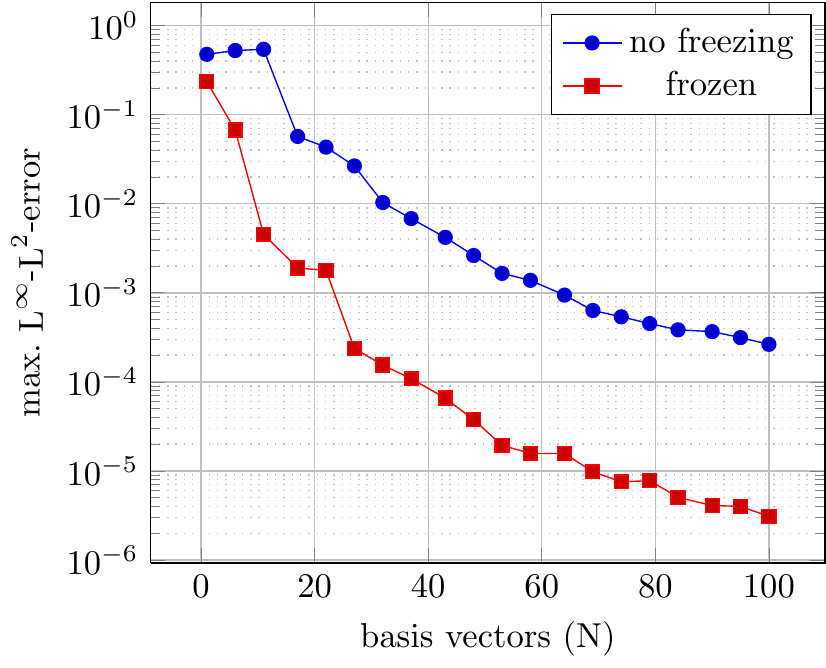} \end{center}
\caption{Error of the FrozenRB-approximation of (\ref{burgersfrozen}) in comparison with the same reduction method without freezing. 
The number of interpolation points $M$ is given by $M = 1.8\,\cdot\,N$.
The error is estimated as the maximum approximation error over a set of 100 randomly chosen parameters.}
\label{fig:errplot}
\end{figure}

Figure \ref{fig:solution} shows that the FrozenRB-scheme provides already for $20$ basis vectors a good approximation of $v_\mu$ (right column) whereas, for the same amount of basis vectors, the solutions of the same scheme without freezing are strongly deformed (center column).
In Figure \ref{fig:errplot} we compare the convergence of the schemes to the corresponding detailed approximations: the FrozenRB-scheme improves the approximation error for equal basis sizes by two orders of magnitude.

\end{document}